\documentclass[english,a4paper,10pt]{amsart}
\usepackage{babel,amsmath,amssymb,amsbsy,amsfonts,latexsym,amsthm, amscd,amsxtra}
\usepackage[T1]{fontenc}
\usepackage{graphics,latexsym}

\def\proc#1{\medbreak\noindent{\it #1}\hspace{1ex}\ignorespaces}
\def\procno#1{\medbreak\noindent{\scshape #1}\hspace{1ex}\ignorespaces}
\newenvironment{Def}{\proc{\bf Definition.}}{\medbreak}
\newtheorem{Lem}{Lemma}
\newtheorem{Prop}{Proposition}
\newtheorem{Cor}{Corollary}
\newenvironment{Claim}{{\bf Claim.}}{\medbreak}
\def\ep{\noindent{\hfill $\Box$}}
\newenvironment{Proof}{\proc{\bf Proof.}}{\ep \medbreak}
\newenvironment{mainTeo}{\procno{\bf Main Theorem.}}{\medbreak}
\newenvironment{mainLem}{\procno{\bf Main Lemma.}}{\medbreak}
\newenvironment{Rem}{\proc{\bf Remark.}}{\medbreak}
\def\acks{\subsubsection*{Acknowledgements.}}

\newcommand\beq[1]{ \begin{equation} \label{#1}}
\newcommand{\eeq}{ \end{equation} }
\newcommand{\beqno}{ \begin{equation*}}
\newcommand{\eeqno}{ \end{equation*}}
\newcommand\beqa[1]{ \begin{eqnarray} \label{#1}}
\newcommand{\eeqa}{ \end{eqnarray} }
\newcommand{\beqano}{ \begin{eqnarray*} }
\newcommand{\eeqano}{ \end{eqnarray*} }

\newcommand{\T}{ {\mathbb T}   }

\newcommand{\R}{ {\mathbb R}   }

\renewcommand{\L}{ {\mathbb L}   }

\renewcommand \a {\alpha}
\newcommand \e {\varepsilon }

\renewcommand \b  {\beta}

\renewcommand \d {\delta}

\newcommand \m {\mu}
\newcommand \n {\nu}

\newcommand \f {\varphi}

\newcommand \g {\gamma}

\newcommand \s {\sigma}

\renewcommand \l {\lambda}

\newcommand \cA {{\mathcal A}}
\newcommand \cS {{\mathcal S}}
\newcommand \cD {{\mathcal D}}
\newcommand \cbA {\bar{{\mathcal A}}}
\newcommand \bd {\bar{\d}}
\newcommand \cL {{\mathcal L}}

\newcommand \dpr {\partial}
\newcommand{\ie}{{\it i.e.}}
\newcommand{\eg}{{\it e.g.}}

\newcommand{\tu}{\tilde{u}}
\newcommand{\tx}{\tilde{x}}

\newcommand \su[1]{\frac{1}{#1}}

\newcommand \ov {\overline}

\newcommand\binomiale[2]{
\left(
\begin{array}{c}
{#1}\\
{#2}
\end{array}
\right)}
\title{On the total disconnectedness of the quotient Aubry set}

\author{Alfonso Sorrentino}
\address{Department of Mathematics, Princeton University, Princeton (NJ), 08544-1000 U.S.}
\email{asorrent@math.princeton.edu}

\begin{document}

\begin{abstract}
In this paper we show that the quotient Aubry set, associated to a sufficiently smooth mechanical or symmetrical Lagrangian,
is totally disconnected (\ie, every connected component consists of a single point). This result is optimal, in the sense of the regularity of the Lagrangian, as Mather's counterexamples in \cite{Mather04} show. Moreover, we discuss the relation between this problem and a Morse-Sard type property for (difference of) critical subsolutions of Hamilton-Jacobi equations.
\end{abstract}

\maketitle


\section{Introduction.}

In Mather's studies of the dynamics of Lagrangian systems
and the existence of Arnold diffusion, it turns
out that understanding certain aspects of the {\it Aubry set} and,
in particular, what is called the {\it quotient Aubry set}, may
help in the construction of orbits with interesting behavior.\\
While in the case of twist maps (see for instance \cite{Bangert,Forni-Mather} and references therein)
there is a detailed
structure theory for these sets, in more degrees of freedom quite few is known. In particular, it seems to
be useful to know whether the quotient Aubry set is ``small'' in
some sense of
dimension (\eg, vanishing topological or box dimension).\\
In \cite{Mather03} Mather showed that if the state space has
dimension $\leq \,2$ (in the non-autonomous case) or the Lagrangian is the kinetic energy
associated to a Riemannian metric and the state space has dimension
$\leq \,3$, then the quotient Aubry set is totally disconnected, \ie, every connected component
consists of a single point (in a compact metric space this is equivalent to vanishing
topological dimension). In the autonomous case, with ${\rm dim}\, M \leq 3$, the same argument shows that this quotient is totally disconnected as long as the Aubry set does not intersect the zero section of ${\rm T}M$ (this is the case when the cohomology class is large enough in norm).\\
What happens in higher dimension? Unfortunately, this is generally
not true. In fact, Burago, Ivanov and Kleiner in \cite{Burago}
provided an example that does not satisfy this property (they did
not discuss it in their work, but it follows from the results
therein). More strikingly, Mather provided in \cite{Mather04}
several examples of quotient Aubry sets that are not only non-totally-disconnected, 
but even isometric to closed intervals. All these
examples come from mechanical Lagrangians on ${\rm T}\T^d$ (\ie, the sum of
the kinetic energy and a potential) with $d\geq 3$. In particular,
for every $\e>0$, he provided a potential $U\in C^{2d-3,1-\e}(\T^d)$,
whose associated quotient Aubry set is isometric to an interval. As
the author himself noticed, it is not possible to improve the
differentiability
of these examples, due to the construction carried out.

The main aim of this article is to show that the counterexamples
provided by Mather are optimal, in the sense that for more regular
mechanical Lagrangians, the associated quotient Aubry sets - corresponding to the zero
cohomology class - are totally disconnected.

In particular, our result will also apply to slightly more general Lagrangians, satisfying certain conditions 
on the zero section; in this case, we shall be able to show that the quotient Aubry set, corresponding to a well specified
cohomology class, is totally disconnected. \\
We shall also outline a possible approach to generalize this result, pointing out how it is related to a {\it Morse-Sard type} problem; from this and Sard's lemma, one can easily recover Mather's result in dimension $d=2$ (autonomous case).

It is important to point out, that most of this
approach has been inspired by Albert Fathi's talk \cite{Fathisem}, in which he used this relation with Sard's lemma to 
show a simpler way to construct mechanical Lagrangians on ${\rm T}\T^N$, whose
quotient Aubry sets are Lipschitz equivalent to any given {\it doubling} metric
space or, equivalently, to any space with finite {\it Assouad dimension} (see \cite{Piotr} for a similar construction). 
In this case we do not get a neat relation between their regularity and $N$, as in
Mather's, but we can only observe that $N$ goes to infinity as $r$
increases. It would be interesting to study in depth the relation
between the dimension of the quotient Aubry set,
the regularity of the Lagrangian and the dimension of the state
space. Our result may be seen as a first step in this direction.

\proc{Post Scriptum.} Just before submitting this paper, we learnt that analogous results had been proven indipendently by Albert Fathi, Alessio Figalli and Ludovic Rifford, using a similar approach (to be published).\\
Moreover, in ``{\it A generic property of families of Lagrangian systems}'' (to appear on {\it Annals of Mathematics}), Patrick Bernard and Gonzalo Contreras managed to show that generically, in Ma\~n\'e's sense, there are at most $1+\dim { H}^1(M;\R)$ ergodic minimizing measures, for each cohomology class $c\in {H}^1(M;\R)$. As a corollary of this striking result, one gets that generically the quotient Aubry set is finite for each cohomology class and it consists of at most $1+\dim {H}^1(M;\R)$ elements.\\


\section{The Aubry set and the quotient Aubry set.}

Let $M$ be a compact and connected smooth manifold without boundary.
Denote by ${\rm T}M$ its tangent bundle and ${\rm T}^*M$ the cotangent one. A
point of ${\rm T}M$ will be denoted by $(x,v)$, where $x\in M$ and $v\in
{\rm T}_xM$, and a point of ${\rm T}^*M$ by $(x,p)$, where $p\in {\rm T}_x^*M$ is a
linear form on the vector space ${\rm T}_xM$. Let us fix a Riemannian
metric $g$ on it and denote with $d$ the induced metric on $M$; let
$\|\cdot\|_x$ be the norm induced by $g$ on ${\rm T}_xM$; we shall use the same notation for the norm induced on
${\rm T}_x^*M$.

\begin{Def}
A function $L:\,{\rm T}M\, \longrightarrow \,\R$ is called a {\it Tonelli
Lagrangian} if:
\begin{itemize}
\item[i)]   $L\in C^2({\rm T}M)$;
\item[ii)]  $L$ is strictly convex in the fibers, \ie, the second partial vertical derivative
$\frac{\dpr^2 L}{\dpr v^2}(x,v)$ is positive definite, as a quadratic form, for any $(x,v)\in {\rm T}M$;
\item[iii)] $L$ is superlinear in each fiber, \ie,
            $$\lim_{\|v\|_x\rightarrow +\infty} \frac{L(x,v)}{\|v\|_x} = + \infty$$
            (this condition is independent of the choice of the Riemannian metric).
\end{itemize}
\end{Def}

Given a Lagrangian, we can define the associated {\it Hamiltonian},
as a function on the cotangent bundle:
\beqano H:\; {\rm T}^*M &\longrightarrow & \R \\
(x,p) &\longmapsto & \sup_{v\in {\rm T}_xM} \{\langle p,\,v \rangle_x -
L(x,v)\}\, \eeqano where $\langle \,\cdot,\,\cdot\, \rangle_x$
represents the canonical pairing between the tangent and cotangent
space.

If $L$ is a Tonelli Lagrangian, one can easily prove that $H$ is
finite everywhere, $C^2$, superlinear and strictly
convex in the fibers. Moreover, under the above assumptions, one can define a
diffeomorphism between ${\rm T}M$ and ${\rm T}^*M$, called the {\it Legendre
transform}:
\beqano \cL:\; {\rm T}M &\longrightarrow & {\rm T}^*M \\
(x,v) &\longmapsto & \left(x,\,\frac{\dpr L }{\dpr v}(x,v)
\right). \eeqano
In particular, $\cL$ is a conjugation between the
two flows (namely the Euler-Lagrange and Hamiltonian flows) and
$$ H \circ \cL(x,v) = \left\langle \frac{\dpr L }{\dpr v}(x,v),\,v \right\rangle_x - L(x,v)\,.$$


Observe that if $\eta$ is a $1$-form on $M$, then we can define a function on the tangent space
\beqano 
\hat{\eta}: {\rm T}M &\longrightarrow&  \R \\
(x,v) &\longmapsto& \langle \eta(x),\, v\rangle_x
\eeqano
and consider a new Tonelli Lagrangian $L_{\eta}= L - \hat{\eta}$. The associated Hamiltonian will be 
$H_{\eta}(x,p) = H(x,p + \eta)$.
Moreover, if $\eta$ is closed, then $\int L\, dt$ and $\int L_{\eta} dt$ have the same extremals and therefore the Euler-Lagrange flows on ${\rm T} M$ associated to $L$ and $L_{\eta}$ are the same.

Although the extremals are the same, this is not generally true for the minimizers. What one can say is that they stay the same when we change the Lagrangian by an exact $1$-form. Thus, for fixed $L$, the minimizers depend only on the de Rham cohomology class $c=[\eta] \in H^1(M;\R)$. 
From here, the interest in considering modified Lagrangians, corresponding to different cohomology classes.\\


Let us fix $\eta$, a smooth ($C^2$ is enough for what follows) $1$-form on $M$, and let $c=[\eta]\in H^1(M;\R)$ be its cohomology class.

As done by Mather in \cite{Mather93}, it is convenient to
introduce, for $t>0$ and $x,\,y \in M$, the following quantity:
$$ h_{\eta,t}(x,y) = \inf \int_0^t L_{\eta}(\g(s),\dot{\g}(s))\,ds\,,$$
where the infimum is taken over all piecewise $C^1$ paths $\g:
[0,t]\longrightarrow M$, such that $\g(0)=x$ and $\g(t)=y$.
We define the {\it Peierls barrier} as:
$$ h_{\eta}(x,y) = \liminf _{t \rightarrow +\infty} (h_{\eta,t}(x,y) + \a(c)t)\,,$$
where $\a:\,H^1(M;\R) \longrightarrow \R$ is Mather's $\a$ function (see \cite{Mather91}). It can be shown that this function is convex and that (only for the autonomous case) the $\liminf$ can be replaced by $\lim$.

Observe that $h_{\eta}$ does not depend only on the cohomology class $c$, but also on the choice of the representant; namely, if $\eta' = \eta + df$, then $h_{\eta'}(x,y) = h_{\eta}(x,y) + f(y) - f(x)$.

\begin{Prop}
The values of the map $h_{\eta}$ are finite. Moreover, the following
properties hold:
\begin{itemize}
\item[{\rm i)}]   $h_{\eta}$ is Lipschitz;
\item[{\rm ii)}]  for each $x\in M$, $h_{\eta}(x,x)\geq 0$;
\item[{\rm iii)}]  for each $x,\,y,\,z \in M$, $h_{\eta}(x,y) \leq h_{\eta}(x,z) + h_{\eta}(z,y)$;
\item[{\rm iv)}]  for each $x,\,y \in M$, $h_{\eta}(x,y) + h_{\eta}(y,x) \geq 0$.
\end{itemize}
\end{Prop}

For a proof of the above claims and more, see \cite{Mather93, Fathibook, ContrerasIturriaga}.
Inspired by these properties, we can define
\beqano \d_c:\; M\times M &\longrightarrow & \R \\
(x,y) &\longmapsto & h_{\eta}(x,y) + h_{\eta}(y,x) \eeqano

(observe that this function does actually depend only on the cohomology class).

This function is positive,
symmetric and satisfies the triangle inequality; therefore, it is a
pseudometric on $$\cA_{L,c} =\{x\in M : \; \d_c(x,x)=0 \}\,. $$

$\cA_{L,c}$ is called the {\it Aubry set} (or {\it projected Aubry set})
associated to $L$ and $c$, and $\d_c$ is {\it Mather's pseudometric}. In
\cite{Mather93}, Mather has showed that this is a non-empty compact
subset of $M$, that can be Lipschitzly lifted to a compact invariant subset of ${\rm
T}M$.

\begin{Def}
The {\it quotient Aubry set} $(\cbA_{L,c},\, \bd_c)$ is the metric space
obtained by identifying two points in $\cA_{L,c}$, if their $\d_c$-pseudodistance is
zero.
\end{Def}

We shall denote an element of this quotient by $\bar{x} = \{y\in
\cA_{L,c}:\; \d_c(x,y)=0 \}$. These elements (that are also called {\it
$c$-static classes}, see \cite{ContrerasIturriaga}) provide a partition of $\cA_{L,c}$ into compact subsets,
that can be lifted to invariant subsets of ${\rm T}M$. They are
really interesting from a dynamical systems point of view, since
they contain the $\a$ and $\omega$ limit sets of $c$-minimizing
orbits (see \cite{Mather93,ContrerasIturriaga} for more details).\\

For the sake of our proof, it is convenient to adopt Fathi's {\it
weak KAM theory} point of view (we remand the reader to \cite{Fathibook} for a self-contained presentation).

\begin{Def}
A locally lipschitz function $u: M \longrightarrow \R$ is a {\it
subsolution} of $H_{\eta}(x,d_xu)=k$, with $k\in \R$, if $H_{\eta}(x,d_xu)\leq k$
for almost every $x\in M$.
\end{Def}

This definition makes sense, because, by Rademacher's theorem, we
know that $d_xu$ exists almost
everywhere.\\
It is possible to show that there exists $c[\eta]\in \R$, such that
$H_{\eta}(x,d_xu)=k$ admits no subsolutions for $k<c[\eta]$ and has
subsolutions for $k\geq c[\eta]$. The constant $c[\eta]$ is called {\it
Ma\~n\'e's critical value} and coincides with $\a(c)$, where $c=[\eta]$ (see
\cite{ContrerasIturriaga}).

\begin{Def}
$u: M \longrightarrow \R$ is a $\eta$-{\it critical subsolution}, if
$H_{\eta}(x,d_xu)\leq \a(c)$ for almost every $x\in M$.
\end{Def}

Denote by $\cS_{\eta}$ the set of critical subsolutions.
This set $\cS_{\eta}$ is non-empty. In fact, Fathi showed (see \cite{Fathibook}) that:

\begin{Prop}
If $u: M \longrightarrow \R$ is a $\eta$-{\it critical subsolution}, then
for every $x,\,y \in M$
$$u(y) - u(x) \leq h_{\eta}(x,y)\,. $$
Moreover, for any $x\in M$, the function $h_{\eta,x} (\cdot) := h_{\eta}(x,\,\cdot)$ is a
$\eta$-critical subsolution.
\end{Prop}

Using this result, he provided a nice representation of $h_{\eta}$, in
terms of the $\eta$-critical subsolutions.

\begin{Cor}
If $x\in \cA_{L,c}$ and $y\in M$,
$$ h_{\eta}(x,y) = \sup_{u \in \cS_{\eta}} (u(y)-u(x))\,.$$
This supremum is actually attained.
\end{Cor}

\begin{Proof}
It is clear, from the proposition above, that
$$h_{\eta}(x,y) \geq \sup_{u\in \cS_{\eta}} (u(y)-u(x))\,.$$
Let us show the other inequality. In fact,
since $h_{\eta,x}$ is a $\eta$-critical subsolution and $x\in \cA_{L,c}$ (\ie,
$h_{\eta}(x,x)=0$), then
$$h_{\eta}(x,y) = h_{\eta,x}(y)-h_{\eta,x}(x) \leq \sup_{u \in \cS_{\eta}} (u(y)-u(x))\,.$$
This shows that the supremum is attained.
\end{Proof}

This result can be still improved. Fathi and Siconolfi proved in \cite{FathiSiconolfi}:\\

{\bf Theorem ({Fathi, Siconolfi}).}
For any $\eta$-critical subsolution $u: M \longrightarrow \R$  and for each
$\e>0$, there exists a $C^1$ function $\tu: M \longrightarrow \R$
such that:
\begin{itemize}
\item[{\rm i)}]   $\tu(x)=u(x)$ and $H_{\eta}(x,d_x\tu)=\a(c)$ on $\cA_{L,c}$;
\item[{\rm ii)}]  $|\tu(x)-u(x)|<\e$ and $H_{\eta}(x,d_x\tu)<\a(c)$ on $M\setminus \cA_{L,c}$.
\end{itemize}
\medbreak

In particular, this implies that $C^1$ $\eta$-critical subsolutions are
dense in $\cS_{\eta}$ with the uniform topology.
This result has been recently improved by Patrick Bernard (see \cite{Bernard}),
showing that every $\eta$-critical subsolution coincides, on the Aubry set, with a 
$C^{1,1}$ $\eta$-critical subsolution.

Denote the set of $C^1$  $\eta$-critical subsolutions by $\cS^1_{\eta}$ and the set of 
$C^{1,1}$  $\eta$-critical subsolutions by $\cS^{1,1}_{\eta}$.

\begin{Cor}
For $x,\,y \in \cA_{L,c}$, the following representation holds:
$$ h_{\eta}(x,y) = \sup_{u \in \cS^{1}_{\eta}} (u(y)-u(x)) = \sup_{u \in \cS^{1,1}_{\eta}} (u(y)-u(x))\,.$$
Moreover, these suprema are attained.
\end{Cor}

It turns out to be convenient, to characterize the elements of $\cbA_{L,c}$ (\ie, the $c$-quotient classes) in terms of
$\eta$-critical subsolutions.

Let us consider the following set:
$$\cD_{c}=\{ u-v:\; u,\,v \in \cS_{\eta}\}\, $$
(it depends only on the cohomology class $c$ and not on $\eta$) and denote by 
$\cD_{c}^{1}$ and $\cD_{c}^{1,1}$, the sets corresponding, respectively, to $C^1$ and $C^{1,1}$ 
$\eta$-critical subsolutions.

\begin{Prop}
For $x,\,y \in \cA_{L,c}$, 
\beqano
\d_c(x,y) &=& \sup_{w \in \cD_{c}} (w(y)-w(x)) = \sup_{w \in \cD_{c}^1} (w(y)-w(x)) =\\
&=& \sup_{w \in \cD_{c}^{1,1}} (w(y)-w(x))\,
\eeqano
and this suprema are attained.
\end{Prop}

\begin{Proof}
From the definition of $\d_c(x,y)$, we immediately get: \beqano
\d_c(x,y) &=& h_{\eta}(x,y) + h_{\eta}(y,x) = \\
&=& \sup_{u \in \cS_{\eta}} (u(y)-u(x)) +
\sup_{v \in \cS_{\eta}} (v(x)-v(y)) = \\
&=& \sup_{u, v \in \cS_{\eta}} [(u(y)-v(y)) - (u(x)-v(x))]  =
\\
&=& \sup_{w \in \cD_{c}} (w(y)-w(x))\,. \eeqano 
The other equalities follow from the density results we mentioned above.
\end{Proof}

\begin{Prop}
If $w \in \cD_{c}$, then $d_xw =0$ on $\cA_{L,c}$. Therefore $ \cA_{L,c} \subseteq
\bigcap_{w\in\cD_{c}^{1,1}} {\rm Crit}(w)\,,$ where ${\rm Crit}(w)$ is the
set of critical points of $w$.
\end{Prop}

\begin{Proof}
This is an immediate consequence of a result by Fathi (see
\cite{Fathibook}); namely, if $u,\,v \in \cS_{\eta}$, then they are
differentiable on $\cA_{L,c}$ and $d_xu = d_xv$.
\end{Proof}

\begin{Prop}
If $w\in \cD_{c}$, then it is constant on any quotient class of
$\cbA_{L,c}$; namely, if $x,\,y \in \cA_{L,c}$ and $\d_c(x,y)=0$, then
$w(x)=w(y)$.
\end{Prop}

\begin{Proof}
From the representation formula above, it follows that: 
\beqano 0
&=& \d_c(x,y) = \sup_{\tilde{w} \in \cD_{c}} (\tilde{w}(y)-\tilde{w}(x))
\geq w(y) -w(x) \\
0 &=& \d_c(y,x) =  \sup_{\tilde{w} \in \cD_{c}}
(\tilde{w}(x)-\tilde{w}(y)) \geq w(x) -w(y)\,. \eeqano
\end{Proof}

For any $w\in \cD_{c}^{1}$, let us define the following {\it evaluation
function}: \beqano \f_w : (\cbA_{L,c},\,\bd_c) &\longrightarrow&
(\R,\,|\cdot|)\\
\bar{x} &\longmapsto& w(x)\,. \eeqano

\begin{itemize}
\item $\f_w$ is well defined, \ie, it does not depend on the element
of the class at which $w$ is evaluated;
\item $\f_w(\cbA_{L,c}) = w(\cA_{L,c}) \subseteq w({\rm Crit}(w))$;
\item $\f_w$ is Lipschitz, with Lipschitz constant $1$.
In fact: \beqano \f_w(\bar{x})-\f_w(\bar{y}) &=& w(x)-w(y) \leq \d_c(x,y) = \bd_c(\bar{x},\bar{y}) \\
\f_w(\bar{y})-\f_w(\bar{x}) &=& w(y)-w(x) \leq \d_c(y,x)
= \bd_c(\bar{y},\bar{x})\,. \eeqano Therefore:
$$|\f_w(\bar{x})-\f_w(\bar{y})| \leq \bd_c(\bar{x},\bar{y})\,.$$
\end{itemize}

As we shall see, these functions play a key role in the proof of our result.

\section{The main result.}

Our main goal is to show that, under suitable hypotheses on $L$, there is a well specified
cohomology class $c_L$, for which $(\cbA_{L,c_L},\,\bd_{c_L})$ is totally disconnected, \ie, every
connected component consists of a single point.

Consider $L: TM \longrightarrow \R$ a Tonelli Lagrangian and the associated Legendre transform
\beqano \cL:\; {\rm T}M &\longrightarrow & {\rm T}^*M \\
(x,v) &\longmapsto & \left(x,\,\frac{\dpr L }{\dpr v}(x,v)
\right). \eeqano

Remember that ${\rm T}^*M$, as a cotangent bundle, may be equipped with a {\it canonical} symplectic structure.
Namely, if $({\mathcal U},\,x_1,\ldots, x_d)$ is a local coordinate chart for $M$ and $({\rm T}^*{\mathcal U},\,x_1,\ldots, x_d, p_1,\ldots, p_d)$ the associated cotangent coordinates, one can define the $2$-form
$$ \omega = \sum_{i=1}^d dx_i \wedge dp_i \,.$$
It is easy to show that $\omega$ is a symplectic form (\ie, it is non-degenerate and closed). In particular, one can check that $\omega$ is intrinsically defined, by considering the $1$-form on ${\rm T}^*{\mathcal U}$
$$ \lambda = \sum_{i=1}^d p_i\,dx_i\,,$$
which satisfies $\omega = -d\lambda$ and is coordinate-indipendent; in fact, in terms of the natural projection 
\beqano 
\pi: {\rm T}^*M &\longrightarrow& M \\
(x,p) &\longmapsto& x
\eeqano
the form $\lambda$ may be equivalently defined pointwise without coordinates by
$$ \lambda_{(x,p)} = (d\pi_{(x,p)})^*p \;\in {\rm T}^*_{(x,p)}{\rm T}^*M\,.$$
The $1$-form $\lambda$ is called the {\it Liouville form} (or the {\it tautological form}).\\

Consider now the section of ${\rm T}^*M$ given by
$$ \Lambda_{L} = \cL(M \times \{0\}) = \left\{ \left(x, \frac{\dpr L }{\dpr v}(x,0) \right):\quad x\in M\right\}\,,$$
corresponding to the $1$-form
$$ \eta_L(x) = \frac{\dpr L }{\dpr v}(x,0) \cdot dx  = \sum_{i=1}^d \frac{\dpr L }{\dpr v_i}(x,0)\,dx_i\,.$$

We would like this $1$-form to be closed, that is equivalent to ask $\Lambda_L$ to be a Lagrangian submanifold, in order to  consider its cohomology class $c_L = [\eta_L ] \in H^1(M;\R)$. Observe that this cohomology class can be defined in a more intrinsic way; in fact, consider the projection 
$$\pi_{|\Lambda_L}: \Lambda_L \subset {\rm T}^*M \longrightarrow M\,;$$ 
this induces an isomorphism between the cohomology groups
$H^1(M;\R)$ and $H^1(\Lambda_L;\R)$. The preimage of $[\lambda_{|\Lambda_L}]$ under this isomorphism is called the
{\it Liouville class} of $\Lambda_L$ and one can easily show that it coincides with $c_L$.\\

We can define the set:
$$ \L(M) = \left\{L:{\rm T}M \longrightarrow \R:\; L\; \text{is a Tonelli Lagrangian and } \Lambda_L\; \text{is Lagrangian}\right\}\,.$$
This set is non-empty and consists of Lagrangians of the form 
$$L(x,v) = f(x) + \langle \eta(x),\,v\rangle_x + O(\|v\|^2), $$
with $f\in C^2(M)$ and $\eta$ a $C^2$ closed $1$-form on $M$.
In particular, it includes the
{\rm mechanical Lagrangians}, \ie, Lagrangians of the form 
$$L(x,v)=\su{2} \|v\|_x^2 + U(x)\,,$$ 
namely the sum of the {\rm kinetic energy} and a {potential} $U: M \longrightarrow \R$. More generally, 
it contains the {\it symmetrical} (or {\it reversible}) {\it Lagrangians}, \ie, Lagrangians $L:\,{\rm T}M\, \longrightarrow \,\R$ such that 
$$ L(x,v)= L(x,-v)\,,$$
for every $(x,v)\in {\rm T}M$.

In fact, in the above cases, $\frac{\dpr L }{\dpr v}(x,0) \equiv 0$; therefore $\Lambda_L = M\times \{0\}$ (the zero section of the cotangent space), that is clearly Lagrangian, and $c_L = 0$.\\

We can now state our main result:

\begin{mainTeo}{\it
Let $M$ be a compact connected manifold of dimension $d\geq 1$ and
let $L \in \L(M)$ be a Lagrangian such that $L(x,0)\in C^r(M)$, with $r\geq 2d-2$
and $\frac{\dpr L}{\dpr v}(x,0)\in C^2(M)$.
Then, the quotient Aubry set $(\cbA_{L,c_L},\,\bd_{c_L})$, corresponding to the {\it Liouville class} of $\Lambda_L$, is totally disconnected, \ie, every connected component consists of a single point.}
\end{mainTeo}

This result immediately implies:

\begin{Cor}[Symmetrical Lagrangians]
Let $M$ be a compact connected manifold of dimension $d\geq 1$ and
let $L(x,v)$ be a symmetrical Tonelli Lagrangian on
${\rm T}M$, such that $L(x,0)\in C^r(M)$, with $r\geq 2d-2$.
Then, the quotient Aubry set $(\cbA_{L,0},\,\bd_0)$ is totally disconnected.
\end{Cor}

More specifically,

\begin{Cor}[Mechanical Lagrangians]
Let $M$ be a compact connected manifold of dimension $d\geq 1$ and
let $L(x,v)=\su{2} \|v\|_x^2 + U(x)$ be a mechanical Lagrangian on
${\rm T}M$, such that the potential $U \in C^r(M)$, with $r\geq 2d-2$.
Then, the quotient Aubry set $(\cbA_{L,0},\,\bd_0)$ is totally disconnected.
\end{Cor}

\begin{Rem}
This result is optimal, in the sense of the regularity of
the potential $U$, for $\cbA_{L,0}$ to be totally disconnected. In fact,
Mather provided in \cite{Mather04} examples of quotient Aubry sets
isometric to the unit interval, corresponding to mechanical Lagrangians $L\in
C^{2d-3,1-\e}({\rm T}\T^d)$, for any $0<\e<1$.
\end{Rem}

Before proving the main theorem, it will be useful to show some useful results.

\begin{Lem}
Let us consider $L\in \L(M)$, such that
$\frac{\dpr L}{\dpr v}(x,0)\in C^2(M)$,
and let $H$ be the associated Hamiltonian. 
\begin{enumerate}
\item Every constant function $u\equiv const$ is a $\eta_L$-critical subsolution. In particular, all $\eta_L$-critical subsolutions are such that
$d_xu\equiv 0$ on $\cA_{L,c_L}$.
\item For every $x\in M$, 
      $$\frac{\dpr H_{\eta_L}}{\dpr p}(x,0) = \frac{\dpr H}{\dpr p}(x,\eta_L(x)) = 0\,.$$
\end{enumerate}
\end{Lem}

\begin{Proof}\
\begin{enumerate}
\item The second part follows immediatly from the fact that, if $u,\,v \in \cS_{\eta_L}$, then they are
differentiable on $\cA_{L,c_L}$ and $d_xu = d_xv$ (see \cite{Fathibook}).\\
Let us show that $u\equiv const$ is a $\eta_L$-critical subsolution; namely, that
$$ H_{\eta_L}(x,0) \leq \a(c_L)$$
for every $x\in M$.
It is sufficient to observe:
\begin{itemize}
\item	$H_{\eta_L}(x,0)= - L(x,0)$; in fact:
	\beqano
	H_{\eta_L}(x,0) &=& H(x,\eta_L(x)) = H\left(x,\frac{\dpr L }{\dpr v}(x,0)\right) = \\
	&=& \left\langle \frac{\dpr L }{\dpr v}(x,0),\,0 \right\rangle_x - L(x,0) = \\
	&=& - L(x,0)\,.
	\eeqano
\item 	let $v$ be {\it dominated} by $L_{\eta_L} + \a(c_L)$ (see \cite{Fathibook}, for the existence of such functions), \ie, for each continuous piecewise $C^1$	curve $\g: [a,b] \longrightarrow M$ we have
	$$ v(\g(b)) - v(\g(a)) \leq \int_a^b L_{\eta_L}(\g(t),\dot{\g}(t))\,dt + \a(c_L)(b-a)\,.$$
	Then, considering the constant path $\g(t)\equiv x$, one can easily deduce that
	\beqano
	\a(c_L) \geq \sup_{x\in M} (-L_{\eta_L}(x,0)) = - \inf_{x\in M} L_{\eta_L}(x,0)\,; 
	\eeqano
	therefore, $$\a(c_L) \geq -L_{\eta_L}(x,0) = - L(x,0) = H_{\eta_L}(x,0)$$ for every $x\in M$.
\end{itemize}
\item	The inverse of the Legendre transform can be written in coordinates
	\beqano \cL^{-1}:\; {\rm T}^*M &\longrightarrow & {\rm T}M \\
	(x,p) &\longmapsto & \left(x,\,\frac{\dpr H }{\dpr p}(x,p)
\right)\,. \eeqano
	Therefore,
	\beqano
	(x,0) &=& \cL^{-1}\left(\cL (x,0\right)) = \cL^{-1}\left(x,\frac{\dpr L }{\dpr v}(x,0)\right) = \\
	&=& \cL^{-1}((x,\eta_L(x))) = \left(x,\,\frac{\dpr H }{\dpr p}(x,\eta_L(x)\right)\,.
	\eeqano
\end{enumerate}
\end{Proof}

In particular, observing that for any $\eta_L$-critical subsolution $u$, 
$H_{\eta_L}(x,d_xu) = \a(c_L)$ on $\cA_{L,c_L}$, we can easily deduce from above that:
$$\cA_{L,c_L} \subseteq \{L(x,0) = -\a(c_L)\} = \{H(x,\eta_L(x)) = \a(c_L)\}\,$$
and
$$ \a(c_L) = \sup_{x\in M} (-L(x,0)) = - \inf_{x\in M} L(x,0)=:e_0\,,$$
as denoted in \cite{ManeI, ManeII}.

Let us observe that in general $$e_0 \leq \min_{c\in H^1(M;\R)} \a(c) = -\b(0)\,,$$
where $\b: H_1(M;\R) \longrightarrow \R$ is Mather's $\b$-function, \ie, the convex conjugate of $\a$ 
(in \cite{ManeI, ManeII}, the right-hand-side quantity is referred to as {\it strict critical value}). 
Therefore, we are considering an extremal case in which 
$e_0 = \a(c_L) = \min \a(c)$; it follows also quite easily that $c_L \in \dpr \b(0)$, namely, 
it is a subgradient of $\b$ at $0$.\\

A crucial step in the proof of our result will be the following lemma, that can be read
as a sort of relaxed version of Sard's Lemma (the proof will be mainly based on
the one in \cite{Abraham}).

\begin{mainLem}\label{mainlemma}{\it
Let $U\in C^{r}(M)$, with $r\geq 2d-2$, be a non-negative function, vanishing
somewhere and denote $\cA=\{U(x)=0\}$. If $u: M \longrightarrow \R$ is $C^1$ and
satisfies $\|d_xu\|_x^2 \leq U(x)$ in an open neighborhood of $\cA$,
then $|u(\cA)|=0$ {\rm(}where $|\cdot|$ denotes the Lebesgue measure in
$\R${\rm)}.}
\end{mainLem}

See section \ref{appenmainlemma} for its proof. 

In particular, it implies this essential property.

\begin{Cor}
Under the hypotheses of the main theorem, if $u\in \cS_{\eta_L}$, then
$$|u(\cA_{L,c_L})|=0$$
{\rm(}where $|\cdot|$ denotes the Lebesgue measure in $\R${\rm)}.
\end{Cor}

{\bf Proof} ({Corollary}).
First of all, we can assume that $u\in \cS_{\eta_L}^{1}$, because of Fathi and Siconolfi's theorem.
By Taylor's formula, it follows that there exists an open neighborhood $W$ of $\cA_{L,c_L}$, such that for all $x\in W$:
\beqano
\a(c_L) &\geq& H_{\eta_L}(x,d_xu) = H_{\eta_L}(x,0) + \frac{\dpr H_{\eta_L}}{\dpr p}(x,0) \cdot d_xu + \\
&\qquad& + \; \int_0^1 (1-t) \frac{\dpr^2 H_{\eta_L}}{\dpr p^2}(x,t\,d_xu)(d_xu)^2\,dt\,.
\eeqano

Let us observe the following.
\begin{itemize}
\item 	From the previous lemma, one has that $$\frac{\dpr H_{\eta_L}}{\dpr p}(x,0)=0\,,$$
	for every $x\in M$. 
\item	From the strict convexity hypothesis, it follows that there exists ${\g}>0$ such that:
	$$ \frac{\dpr^2 H}{\dpr p^2}(x,t\,d_x u)(d_x u)^2 \geq \g \|d_xu\|_x^2 $$
	for all $x\in M$ and $0\leq t \leq 1$.
\end{itemize}

Therefore, for $x\in W$:
\beqano
\a(c_L) &\geq& H_{\eta_L}(x,d_xu) \geq 
H_{\eta_L}(x,0) + \frac{\g}{2}\|d_xu\|_x^2 =\\
&=& -L(x,0) + \frac{\g}{2}\|d_xu\|_x^2 \,.
\eeqano

The assertion will follow from the previous lemma, choosing 
$$U(x)= \frac{2}{\g}(\a(c_L)+L(x,0))\,.$$
In fact, $U\in C^{r}$, with $r\geq 2d-2$, by hypothesis; moreover, it satisfies all other conditions, because
$$\a(c_L)= - \inf_{x\in M} L(x,0)$$ and 
$$\cA_{L,c_L}\subseteq\{x\in W:\; L(x,0)=-\a(c_L)\} = \{x\in W:\; U(x)=0\}=:\cA\,.$$
For, the previous lemma allows us to conclude that
$$|u(\cA_{L,c_L})| =0\,.$$
\ep
\medbreak

{\bf Proof} ({Main Theorem}).
Suppose by contradiction that $\cbA_{L,c_L}$ is not totally disconnected;
therefore it must contain a connected component $\ov{\Gamma}$ with
at least two points $\bar{x}$ and $\bar{y}$. In particular
$$\bd_c(\bar{x},\bar{y})=h_{\eta_L}(x,y) + h_{\eta_L}(y,x) > 0,$$
for some $x\in \bar{x}$ and $y\in \bar{y}$;
therefore, we have $h_{\eta_L}(x,y)>0$ or $h_{\eta_L}(y,x)>0$.
From the representation formula for $h_{\eta_L}$, it follows that there exists $u\in \cS_{\eta_L}^{1,1} \subseteq \cD_{c_L}^{1,1}$ (since $u=u-0$, and $v=0$ is a $\eta_L$-critical subsolution), such that
$|u(y)-u(x)|>0$. 

This implies that the set $\f_u(\ov{\Gamma})$ is a connected set in $\R$
with at least two different points, hence it is a non degenerate
interval and its Lebesgue measure is positive. But
$$ \f_u(\ov{\Gamma}) \subseteq \f_u(\cbA_{L,c_L}) = u(\cA_{L,c_L})$$
and consequently
$$ 0< \left|\f_u(\ov{\Gamma})\right| \leq |u(\cA_{L,c_L})|.$$
This contradicts the previous corollary.
\ep
\medbreak

In particular, this proof suggests a possible approach to generalize the above result to more general Lagrangians and other cohomology classes.

\begin{Def}
A $C^1$ function $f:M \longrightarrow \R$ is of {\it Morse-Sard type} if
$|f({\rm Crit}(f))|=0$, where ${\rm Crit}(f)$ is the set of critical points of $f$ and 
$|\cdot|$ denotes the Lebesgue measure in $\R$.
\end{Def}

\begin{Prop}
Let $M$ be a compact connected manifold of dimension $d\geq 1$, $L$ a Tonelli Lagrangian and $c\in H^1(M;\R)$.
If each $w\in\cD_{c}^{1,1}$ is of {\it Morse-Sard type}, then
the quotient Aubry set $(\cbA_{L,c},\,\bd_{c})$ is totally disconnected.
\end{Prop}

This proposition and Sard's lemma (see \cite{Bates}) easily imply Mather's result in dimension $d\leq 2$ (autonomous case); it suffices to notice that Sard's lemma (in dimension $d$) holds for $C^{d-1,1}$ functions.

\begin{Cor}
Let $M$ be a compact connected manifold of dimension $d\leq 2$. For any $L$ Tonelli Lagrangian and $c\in H^1(M;\R)$,
the quotient Aubry set $(\cbA_{L,c},\,\bd_{c})$ is totally disconnected.
\end{Cor}

\begin{Rem}
The main problem becomes now to understand under which conditions on $L$ and $c$, these differences of subsolutions are of {\it Morse-Sard type}. 
Unfortunately, one cannot use the classical Sard's lemma, due to a lack of regularity of critical subsolutions: in general they will be at most $C^{1,1}$. In fact, although it is always possible to smooth them up out of the Aubry set and obtain functions in $C^{\infty}(M\setminus \cA_{L,c})\cap C^{1,1}(M)$, the presence of the Aubry set (where the value of their differential is prescribed) represents an obstacle that it is impossible to overcome. It is quite easy to construct examples that do not admit $C^2$ critical subsolutions: just consider a case in which $\cA_{L,c}$ is all the manifold and it is not a $C^1$ graph. For instance, this is the case if $M=\T$ and 
$H(x,p)=\frac{1}{2}\left(p+ \frac{2}{\pi}\right)^2 - \sin^2 (\pi x)$;  in fact, there is only one critical subsolution (up to constants), that turns out to be a solution ($\cA_{L,{\frac{2}{\pi}}} = \T$), and it is given by a primitive of $\sin (\pi x) - \frac{2}{\pi}$; this is clearly $C^{1,1}$ but not $C^2$. \\
On the other hands, the above results suggest that, in order to prove the Morse-Sard property, one could try to control the {\it complexity} of these functions ({\it \`a la } Yomdin), using the rigid structure provided by Hamilton-Jacobi equation and the smoothness of the Hamiltonian, rather than the regularity of the subsolutions. There are several difficulties in pursuing this approach in the general case, mostly related to the nature of the Aubry set. We hope to understand these ``speculations'' more in depth in the future.
\end{Rem}


\section{Proof of the Main Lemma.}\label{appenmainlemma}

\begin{Def}
Consider a function $f\in C^r(\R^d)$. We say that f is {\it s -
flat} at $x_0\in \R^d$ (with $s\leq r$), if all its derivatives, up
to the order $s$, vanish at $x_0$.
\end{Def}

The proof of the main lemma is based on the following
version of {\it Kneser-Glaeser's Rough composition theorem}
(see \cite{Abraham, Whitney}).

\begin{Prop}\label{roughcomp}
Let $V,\,W \subset \R^d$ be open sets, $A\subset V$, $A^* \subset W$
closed sets.
Consider $U\in C^r (V)$, with $r\geq 2$, a non-negative function
that is $s$-flat on $A\subset \{U(x)=0\}$, with $s\leq r-1$, and $g:W \longrightarrow
V$ a $C^{r-s}$ function, with $g(A^*)\subset A$.\\
Then, for every open pre-compact set $W_1 \supset A^*$ properly contained in $W$,
there exists $$F:\R^d \longrightarrow \R$$ satisfying the
following properties:
\begin{itemize}
\item[{\rm i)}] $F\in C^{r-1}(\R^d)$;
\item[{\rm ii)}] $F\geq 0$;
\item[{\rm iii)}] $F(x)=U(g(x))=0$ on $A^*$;
\item[{\rm iv)}] $F$ is $s$-flat on $A^*$;
\item[{\rm v)}] $\{F(x)=0\}\cap W_1 = A^* $;
\item[{\rm vi)}] there exists a constant $K>0$, such that $U(g(x))\leq K F(x)$
on $W_1$.
\end{itemize}
\end{Prop}
See section \ref{appendrough} for its proof.\\

To prove the main lemma, it will be enough to show that for every $x_0 \in M$, there exists a neighborhood $\Omega$
such that it holds. For such a local result, we can assume that $M={\mathcal U}$ is an open subset of $\R^d$, with $x_0 \in {\mathcal U}.$ In the sequel, we shall identify ${\rm T}^*{\mathcal U}$ with ${\mathcal U}\times \R^d$ and for $x\in {\mathcal U}$, we identify ${\rm T}_x^*{\mathcal U} = \{x\} \times \R^d$. We equip ${\mathcal U}\times \R^d$ with the natural coordinates $(x_1,\ldots, x_d,p_1,\ldots, p_d)$.\\
Before proceeding in the proof, let us point out that it is locally possible to replace the norm obtained by the
Riemannian metric, by a constant norm on $\R^d$.

\begin{Lem}\label{norme}
For each $ 0<\a<1$ and $x_0 \in M$, there exists an open neighborhood $\Omega$ of $x_0$, with $\ov{\Omega} \subset {\mathcal U}$ and such that
$$ (1-\a) \|p\|_{x_0} \leq \|p\|_x \leq (1+\a) \|p\|_{x_0}\,, $$
for every $p \in {\rm T}_x^*{\mathcal U} \cong \R^d$ and each $x\in \ov{\Omega}.$
\end{Lem}

\begin{Proof}
By continuity of the Riemannian metric, the norm $\|p\|_x$ tends uniformly to $1$ on $\{p:\; \|p\|_{x_0}=1\}$, as $x$ tends to $x_0$. Therefore, for $x$ near to $x_0$ and every $p \in \R^d\setminus\{0\}$, we have:
$$ (1-\a) \leq \left\| \frac{p}{\|p\|_{x_0}}  \right\|_x  \leq (1+\a). $$
\end{Proof}

We can now prove the main result of this section.

{\bf Proof} ({ Main Lemma}).
By choosing local charts and by lemma \ref{norme}, we can assume that $U\in C^{r}(\Omega)$, with $\Omega$ open set in $\R^d$,
$\cA=\{x \in \Omega:\; U(x)=0\}$ and $u: \Omega \longrightarrow \R$ is
such that $\|d_xu\|^2 \leq \b U(x)$ in $\Omega$, where $\b$ is a positive constant.

Define, for $1\leq s \leq r$:
$$B_s = \{x\in \cA:\; U\;{\rm is }\; s\; \text{- flat at }x\}$$
and observe that
$$\cA = B_1 := \{x\in \cA: DU(x)=0\}\,. $$

We shall prove the lemma by induction on the dimension $d$. Let us start with the following claim.

\begin{Claim}
If $s\geq 2d-2$, then $|u(B_s)|=0$.
\end{Claim}

\begin{Proof}
Let $C \subset \Omega$ be a closed cube with edges parallel to the coordinate axes. We shall show that $|u(B_s\cap C)|=0$.
Since $B_s$ can be covered by countably many such cubes, this will prove that $|u(B_s)|=0$.

Let us start observing that, by Taylor's theorem, for any $x\in B_s\cap C$ and $y\in C$ we have $$ U(y)=R_s(x;y),$$ where
$R_s(x;y)$ is Taylor's remainder. Therefore, for any $y\in C$
$$ U(y) = o(\|y-x\|^{s})\,. $$
Let $\l$ be the length of the edge of $C$. Choose an integer $N>0$ and subdivide $C$ in $N^d$ cubes $C_i$ with
edges $\frac{\l}{N}$, and order them so that, for $1\leq i \leq N_0 \leq N^d$, one has $C_i\cap B_s \neq \emptyset$.
Hence,
$$ B_s \cap C = \bigcup_{i=1}^{N_0} B_s \cap C_i.$$
Observe that for every $\e>0$, there exists $\n_0=\n_0(\e)$ such that, if $N\geq \n_0$, $x\in B_s \cap C_i$ and $y\in C_i$, for some $0\leq i \leq N_0$, then
$$ U(y) \leq \frac{\e^2}{4\b (d\l^2)^d} \|y-x\|^s\,.$$

Fix $\e>0$. Choose $x_i \in B_s \cap C_i$ and call $y_i=u(x_i)$. Define, for $N\geq \n_0$, the following intervals in $\R$:
$$ E_i =\left[ y_i - \frac{\e}{2N^d},\; y_i + \frac{\e}{2N^d}\right].$$

Let us show that, if $N$ is sufficiently big, then
$u(B_s \cap C) \subset \bigcup_{i=1}^{N_0} E_i$.\\
In fact, if $x\in B_s\cap C$, then there exists $1\leq i \leq N_0$, such that $x\in B_s\cap C_i$. Therefore,
\beqano
|u(x)-y_i| &=& |u(x)-u(x_i)| = \\
&=&\|d_xu(\tilde{x})\|\cdot\|x-x_i\| \leq \\
&\leq& \sqrt{\b U(\tilde{x})} \|x-x_i\| \leq\\
&\leq& \sqrt{\b\frac{\e^2}{4\b (d\l^2)^d}} \|\tilde{x}-x_i\|^{\frac{s}{2}} \|x-x_i\| \leq \\
&\leq& \frac{\e}{2(d\l^2)^{\frac{d}{2}}}\|x-x_i\|^{\frac{s+2}{2}} \leq\\
&\leq& \frac{\e}{2(d\l^2)^{\frac{d}{2}}} \left(\sqrt{d}\frac{\l}{N}\right)^{\frac{s+2}{2}}\!\!\!\!\!\!\!\!,
\eeqano
where $\tilde{x}$ is a point in the segment joining $x$ and $x_i$.
Since by hypothesis $s\geq 2d-2$, then $\frac{s+2}{2}\geq d$. Hence, assuming that $N>\max\{\l\sqrt{d},\; \n_0\}$,
one gets
$$|u(x)-y_i| \leq \frac{\e}{2N^d}$$
and can deduce the inclusion above.

To prove the claim, it is now enough to observe:
\beqano
|u(B_s \cap C)| &\leq& \left|\bigcup_{i=1}^{N_0} E_i\right| \leq \sum_{i=1}^{N_0} |E_i| \leq\\
&\leq& \e N_0 \frac{1}{N^d}\leq \\
&\leq& \e N^d \frac{1}{N^{d}} =\\
&=& \e\,.
\eeqano
From the arbitrariness of $\e$, the assertion follows easily.
\end{Proof}

This claim immediately implies that $u(B_{2d-2})$ has measure zero.

In particular, this proves the case $d=1$ (since in this case $2d-2=0$) and allows us to start the induction.\\
Suppose to have proven the result for $d-1$ and show it for $d$.
Since
$$ \cA = (B_1\setminus B_2) \cup (B_2\setminus B_3) \cup \ldots \cup (B_{2d-3}\setminus B_{2d-2}) \cup B_{2d-2}\,,$$
it remains to show that $|u(B_s\setminus B_{s+1})|=0$ for $1\leq s \leq 2d-3 \leq r-1$.

\begin{Claim}
Every $\tx \in B_s\setminus B_{s+1}$ has a neighborhood $\tilde{V}$, such that
$$|u((B_s\setminus B_{s+1}) \cap \tilde{V})|=0\,.$$
\end{Claim}

Since $B_s\setminus B_{s+1}$ can be covered by countably many such neighborhoods, this implies
that $u(B_s\setminus B_{s+1})$ has measure zero.

\begin{Proof}
Choose $\tx \in B_s\setminus B_{s+1}$. By definition of these sets, all partial derivatives of order $s$ of $U$ vanish at this point, but there is one of order $s+1$ that does not. Assume (without any loss of generality) that there exists a function
$$ w(x) = \dpr_{i_1} \dpr_{i_2} \ldots \dpr_{i_s} U (x)$$
such that
$$ w(\tx)=0 \qquad {\rm but} \qquad \dpr_1 w(\tx) \neq 0\,.$$
Define
\beqano
h: \Omega &\longrightarrow& \R^d\\
 x &\longmapsto& (w(x),\,x_2,\,\ldots,\,x_d)\,,
\eeqano
where $x=(x_1,\,x_2,\,\ldots,\,x_d)$.
Clearly, $h\in C^{r-s}(\Omega)$ and $Dh(\tx)$ is non-singular; hence, there is an open neighborhood $V$ of $\tx$ such that
$$ h: V \longrightarrow W$$
is a $C^{r-s}$ isomorphism (with $W=h(V)$).\\
Let $V_1$ be an open precompact set, containing $\tilde{x}$ and properly contained in $V$, and define
$A = B_s\cap \ov{V_1}$, $A^*=h(A)$ and $g=h^{-1}$.
If we consider $W_1$, any open set containing $A^*$ and properly contained in $W$,
we can apply proposition \ref{roughcomp} and deduce the existence of
$F:\R^d \longrightarrow \R$ satisfying properties i)-vi).\\
Define
$$\hat{W} = \{(x_2,\,\ldots,\,x_d) \in \R^{d-1}:\; (0,\,x_2,\,\ldots,\,x_d)\in W_1 \} $$
and
$$\hat{U}(x_2,\,\ldots,\,x_d) = C\,F(0,x_2,\ldots,x_d),$$
where $C$ is a positive constant to be chosen sufficiently big. Observe that $\hat{U}\in C^{r-1}(\R^{d-1})$.\\
Moreover, property v) of $F$ and the fact that $A^*=h(A) \subseteq \{0\}\times \hat{W}$
imply that:
$$ A^* = \{0\} \times \hat{B_1}\,,$$
where $\hat{B_1} = \{(x_2,\ldots\, x_d)\in \hat{W}:\; F(0,x_2,\ldots,x_d)=0\}$.
Denote
$$ \hat{\cA} := \{(x_2,\ldots,x_d) \in \hat{W}:\; \hat{U} = 0\}  = \hat{B_1}\,$$
and define the following function on $\hat{W}$:
$$\hat{u} (x_2,\,\ldots,\,x_d) = u(g(0,x_2,\ldots,x_d)).$$
We want to show that these functions satisfy the hypotheses for the $(d-1)$-dimensional case. In fact:
\begin{itemize}
\item $\hat{U}\in C^{r-1}(\R^{d-1})$, with $r-1\geq 2d-3 > 2(d-1)-2$;
\item $\hat{u} \in C^1(\hat{W})$ (since $g$ is in $C^{r-s}(W)$, where $1\leq s\leq r-1$);
\item if we denote by $\m= \sup_{W_1} \|d_x g\| < +\infty$ (since $g$ is $C^1$ on $\ov{W_1}$), then we have that for every point in $\hat{W}$:
\beqano
\|d\hat{u} (x_2,\ldots,x_d)\|^2 &\leq& \|d_x{u} (g(0,x_2,\ldots,x_d))\|^2\| d_xg(0,x_2,\ldots,x_d)\|^2 \leq\\
&\leq& \m^2 \|d_x{u} (g(0,x_2,\ldots,x_d))\|^2 \leq \\
&\leq& \b \m^2 U(g(0,x_2,\ldots,x_d)) \leq \\
&\leq& \b \m^2 K F(0,x_2,\ldots,x_d) \leq\\
&\leq& \hat{U} (x_2,\ldots,x_d),
\eeqano
if we choose $C> \b \m^2K$, where $K$ is the positive constant appearing in proposition \ref{roughcomp}, property vi).
\end{itemize}

Therefore, it follows from the inductive hypothesis, that:
$$ |\hat{u}(\hat{\cA})| = 0.$$
Since,
\beqano
u(B_s\cap V_1) &\subseteq& u(A) = u(g(A^*)) = u(g(\{0\}\times\hat{B_1})) = \\
&=& \hat{u}(\hat{B_1}) = \hat{u}(\hat{\cA})\,, \eeqano
defining $\tilde{V} = V_1$, we may conclude that
$$ |u(B_s\cap \tilde{V})| \leq |\hat{u}(\hat{\cA})| = 0 \,.$$
\end{Proof}
This completes the proof of the Main Lemma.
\ep
\medbreak


\section{Proof of a modified version of Kneser-Glaeser's Rough composition theorem.} \label{appendrough}

Now, let us prove proposition \ref{roughcomp}. We shall mainly follow
the presentation in \cite{Abraham}, adapted to our needs.

{\bf Proof} ({Proposition \ref{roughcomp}}).
Let us start, defining a family of polynomials. Supposing that $g$
is $C^r$ and using the $s$-flatness hypothesis, we have, for $x\in
A^*$ and $k=0,\,1,\ldots,\, r\,$: \beqa{polynwhitney} f_k(x) =
D^k(U\circ g)(x) = \sum_{s < q \leq k}\sum \s_k D^q
U(g(x))\,D^{i_1}g(x) \ldots D^{i_q}g(x)\,, \eeqa
where the second sum is over all the
$q$-tuples of integers $i_1,\,\ldots,\,i_q\, \geq 1$ such that
$i_1+\ldots + i_q = k$, and $\s_k=\s_k(i_1,\ldots,\,i_q)$.\\
The crucial observation is that (\ref{polynwhitney})
makes sense on $A^*$, even when $g$ is $C^{r-s}$ smooth (in fact
$i_j\leq k-q+1 \leq r-s$).

We would like to proceed in the fashion of {\it Whitney's extension theorem}, in order to find a smooth function $F$ such that $D^k F=
f_k$ on $A^*$, and satisfying the stated conditions.

\begin{Rem}
Note that, without any loss of generality,
we can assume that $W$ is contained in an open ball of
diameter $1$. The general case will then follow from this special
one, by a straightforward partition of unity argument.
\end{Rem}

Let us start with some technical lemmata.

\begin{Lem}\label{lemmwhit}
For $x,\,x',\,x_0\in A^*$ and $k=0,\,\ldots,\,r$, we have:
$$ f_k(x')= \sum_{i\leq r-k} \frac{f_{k+i}(x)}{i!}(x'-x)^i +
R_k(x,x')\,,$$ with
$$
\frac{|R_k(x,x')|}{\|x'-x\|^{r-k}} \longrightarrow 0
$$
as $x,\,x' \longrightarrow x_0$ in $A^*$.
\end{Lem}

The proof of this lemma appears without any major modification in \cite{Abraham} (on pages $36$-$37$).\\

Define, for $x\in A^*$ and $y\in \R^d$
$$
P(x,y)=\sum_{i=s+1}^{r} \frac{f_i(x)}{i!}(y-x)^i\,
$$
and its $k$-th derivative
$$
P_k(x,y)=\sum_{i\leq r-k} \frac{f_{i+k}(x)}{i!}(y-x)^i\,.
$$

\begin{Lem}\label{stimaresto}
For $x\in A^*$ and $y\in {W_1}$,
$$ U(g(y)) = P(x,y) + R(x,y)\,,$$
where $|R(x,y)| \leq C \|y-x\|^r$.
\end{Lem}

\begin{Proof}
The proof follows the same idea of lemma \ref{lemmwhit}. By
Taylor's formula for $U$,
$$U(g(y))= \sum_{q=s+1}^r
\frac{D^qU(g(x))}{q!} (g(y)-g(x))^q  + I(g(x),g(y)) (g(x)-g(y))^r\,.
$$
Obviously,
$$\left| I(g(x),g(y)) (g(x)-g(y))^r\right| \leq C_1 \|y-x\|^r\,,$$
therefore it is sufficient to estimate the first term.\\
Observe that:
$$ g(y) = g(x) + \sum_{i=1}^{r-s} D^i g(x)(y-x)^i +
J(x,y)(y-x)^{r-s}\,.
$$
Hence, the first term in the sum above becomes: \beqano &&
\sum_{q=s+1}^r \frac{D^q U(g(x))}{q!} \left[ \sum_{i=1}^{r-s} D^i
g(x)(y-x)^i + J(x,y)(y-x)^{r-s} \right]^q= \\
&& \qquad = \; \sum_{k=s+1}^r a_k
(y-x)^k + \hat{R}(x,y) =\\
&& \qquad = \; P(x,y) + \hat{R}(x,y)\,, \eeqano
 since
$$ a_k = \sum_{s+1\leq q \leq k} \sum D^q U(g(x)) D^{i_1}g(x) \ldots D^{i_q}g(x) = \frac{f_k(x)}{k!}\,.$$
The remainder terms consist of:
\begin{itemize}
\item terms containing $(y-x)^k$, with $k>r$;
\item terms of the binomial product, containing $J(x,y)(y-x)^{r-s}$.
They are of the form:
$$ \ldots\, (y-x)^{(r-s)j + \sum_{i=1}^{r-s} i\a_i}$$
where $\a_i\geq 0$ and $\sum \a_i = q-j$.
Since $q\geq s+1$ and $s\leq r-1$, then:
 \beqano (r-s)j
+ \sum_{i=1}^{r-s} i\a_i &\geq& (r-s)j +
\sum_{i=1}^{r-s} \a_i  = \\
&=& (r-s)j + q-j = \\
&=& rj - sj + q-j \geq\\
&\geq& rj - (s+1)j + s+1 = \\
&=& r + r(j-1) - (s+1)(j-1) =\\
&=& r + (r-s-1)(j-1) \geq r\,. \eeqano
\end{itemize}

Therefore, for $x\in A^*$ and $y\in W_1$
$$ \left|\hat{R}(x,y)\right| \leq C_2 \|y-x\|^r\,,$$
and the lemma follows taking $C=C_1 + C_2$.
\end{Proof}

Next step will consist of creating a {\it Whitney's partition}. We
will start by covering $W_1\setminus A^*$ with an infinite
collection of cubes $K_j$, such that the size of each $K_j$ is roughly
proportional to its distance from $A^*$.

First, let us fix some notation. We shall write $a\prec b$ instead
of ``there exists a positive real constant $M$, such that $a\leq Mb$
''
and $a \approx b$ as short for $a\prec b$ and $b\prec a$.\\
Let $\l=\frac{1}{4\sqrt{d}}$; this choice will come in handy later. For any closed cube $K$ (with edges
parallel to the coordinate axes), $K^{\l}$ will denote the $(1+\l)$
- dilation of $K$ about its center.\\
Let $\|\cdot\|$ be the euclidean metric on $\R^d$ and
$$d(y)=d(y,A^*)= \inf\{\|y-x\|:\; x\in A^*\}\,.$$
If $\{K_j\}_j$ is the sequence of closed
cubes defined below, with edges of length $e_j$, let $d_j$ be its
distance from $A^*$, \ie,
$$d_j=d(A^*,K_j)= \inf \{\|y-x\|:\; x\in A^*,\, y\in K_j\}\,.$$

One can show the following classical lemma (see for instance \cite{Abraham} for a proof).

\begin{Lem}\label{partwhit}
There exists a sequence of closed cubes $\{K_j\}_j$ with edges
parallel to the coordinate axes, that satisfies the following
properties:
\begin{itemize}
\item[{\rm i)}]     the interiors of the $K_j$'s are disjoint;
\item[{\rm ii)}]    $W_1\setminus A^* \subset \bigcup_j K_j$;
\item[{\rm iii)}]   $e_j \approx d_j$;
\item[{\rm iv)}]    $e_j \approx d(y)$ for all $y\in K_j^{\l}$;
\item[{\rm v)}]     $e_j \approx d(z)$ for all $z\in W_1\setminus
A^*$, such that the ball with center $z$ and radius
$\frac{1}{8}d(z)$ intersects $K_j^{\l}$;
\item[{\rm vi)}]    each point of $W_1\setminus A^*$ has a
neighborhood intersecting at most $N$ of the $K_j^{\l}$, where $N$
is an integer depending only on $d$.
\end{itemize}
\end{Lem}

Now, let us construct a partition of unity on $W_1\setminus A^*$. Let
$Q$ be the unit cube centered at the origin. Let $\eta$ be a
$C^{\infty}$ bump function defined on $\R^d$ such that
$$ \eta (y) = \left\{
\begin{array}{ll}
1 & {\rm for}\; y\in Q\\
0 & {\rm for}\; y\not\in Q^{\l}
\end{array} \right.
$$
and $0\leq \eta \leq 1$. Define
$$\eta_j(y)= \eta\left(\frac{y-c_j}{e_j} \right)\,, $$
where $c_j$ is the center of $K_j$ and $e_j$ is the length of its
edge, and consider
$$\s(y)=  \sum_j \eta_j(y)\,.$$
Then, $1\leq \s(y) \leq N $ for all $y\in W_1\setminus A^*$.
Clearly, for each $k=0,\,1,\,2,\,\ldots$ we have that $D^k\eta_j(y)
\prec e_j^{-k}$, for all $y\in W_1\setminus A^*$. Hence, by
properties iv) and vi) of lemma \ref{partwhit}, we have that for
each $k=0,\,1,\,\ldots,\,r$:
$$D^k\eta_j(y) \prec d(y)^{-k} \qquad \text{for all}\; y\in W_1\setminus A^*$$
and
$$D^k\s(y) \prec d(y)^{-k} \qquad \text{for all}\; y\in W_1\setminus A^*\,.$$
Define
$$ \f_j(y) = \frac{\eta_j(y)}{\s(y)}\,.$$
These functions satisfy the following properties:
\begin{itemize}
\item[{\rm i)}] each $\f_j$ is $C^{\infty}$ and supported on $K_j^{\l}$;
\item[{\rm ii)}] $0\leq \f_j(y) \leq 1$ and $\sum_j \f_j(y)=1$, for all $y\in W_1\setminus
A^*$;
\item[{\rm iii)}] every point of $W_1\setminus A^*$ has a
neighborhood on which all but at most $N$ of the $\f_j$'s vanish
identically;
\item[{\rm iv)}] for each $k=0,\,1,\,\ldots,\,r$, $D^k\f_j(y) \prec d(y)^{-k}$ for all $y\in
W_1\setminus A^*$; namely, there are constants $M_k$ such that
$D^k\f_j(y) \leq M_kd(y)^{-k}$;
\item[{\rm v)}] there is a constant $\a$ and points $x_j \in A^*$,
such that:
$$ \|x_j-y\| \leq \a d(y), \qquad {\rm whenever}\; \f_j(y)\neq 0\,.$$
This follows from properties iii) and iv) of lemma \ref{partwhit}.
\end{itemize}

We can now construct our function $F$. Observe that, from lemma
\ref{stimaresto}:
$$ 0\leq U(g(y))= P(x_j,y) + R(x_j,y) \leq P(x_j,y) + C \|y-x_j\|^r\,;$$
therefore $P(x_j,y) \geq -C \|y-x_j\|^r$. \\
First, define
$$ \hat{P}_j(y)= P(x_j,y) + 2C \|y-x_j\|^r$$
where $C$ is the same constant as in lemma \ref{stimaresto}; for
what said above,
\beqa{defPhatj}
\hat{P}_j(y) \geq C\|y-x_j\|^r >0 \qquad {\rm
in}\; W_1\setminus\{x_j\}\,.
\eeqa
Hence, construct $F$
in the following way:
$$F(y) =
\left\{\begin{array}{ll} 0 & y\in A^*\\
\sum_j \f_j(y) \hat{P}_j(y) & y \in \R^d\setminus A^*\,.
\end{array} \right.
$$
We claim that this satisfies all the stated properties i)-vi). In
particular, properties ii), iii) and v) follow immediately from the
definition of $F$ and (\ref{defPhatj}). Moreover, $F\in
C^{\infty}(\R^d\setminus A^*)$. We need to show that $D^kF=f_k$ (for
$k=0,\,1,\,\ldots,\,r-1$) on $\dpr A^*$ (namely, the boundary of $A^*$) and that
$D^{r-1}F$ is continuous on it. The main difficult in the proof, is
that $D^kF$ is expressed as a sum containing terms
$$ D^{k-m}\f_j(y) P_m(x_j,y),$$ where $\f_j(y)\neq 0$. Even if $y$ is
close to some $x_0\in A^*$, it could be closer to $A^*$ and hence
the bound given by property iv) of $\f_j$ might become large. One
can overcome this problem by choosing a point $x^* \in A^*$, so that
$\|x^*-y\|$ is roughly the same as $d(y)$ and hence, $x_j$ is close
to $x^*$.

\begin{Lem}\label{lemma5.2}
For every $\eta>0$, there exists $\d>0$ such that for all $y\in
W_1\setminus A^*$, $x,\,x^*\in A^*$ and $x_0\in \dpr A^*$, we have
$$ \|P_k(x,y) - P_k(x^*,y)\|\leq \eta \,d(y)^{r-k} \leq \eta \|y-x_0\|^{r-k}, $$
whenever $k\leq r$ and
$$
\left\{ \begin{array}{l}
\|y-x\|<\a d(y)\\
\|y-x^*\|<\a d(y)\\
\|y-x_0\|<\d\,,
\end{array}\right.
$$
where $\a$ is the same constant as in {\rm v)} above.
\end{Lem}

See \cite{Abraham} (on page $126$) for its proof.

%
%
%

\begin{Lem}\label{lemma5.3} For every $\eta>0$, there exist $0<\d<1$ and a constant $E$, such that
for all $y\in W_1\setminus A^*$, $x^*\in A^*$ and $x_0\in \dpr A^*$,
we have
$$ \|D^kF(y) - P_k(x^*,y)\|\leq E\, d(y)^{r-k} \leq \eta\, d(y)^{r-k-1},$$
whenever $k\leq r-1$ and
$$
\left\{ \begin{array}{l}
\|y-x^*\|<\a d(y)\\
\|y-x_0\|<\d\,.
\end{array}\right.
$$
\end{Lem}

\begin{Proof}
Let
$$S_{j,k}(x^*,y)= \dpr_k \hat{P}_j(y) - P_k(x^*,y)\,. $$
From lemma \ref{lemma5.2} (with $\eta = \e$, to be defined later)
and the definition of $\hat{P}_j$, we get: \beqano
\|S_{j,k}(x^*,y)\| &\leq& \|\dpr_k \hat{P}_j(y) - P_k(x_j,y)\| +
\|P_k(x_j,y)-P_k(x^*,y)\| \leq \\
&\leq& C_k d(y)^{r-k} + \e d(y)^{r-k} = \\
&=& (C_k  + \e ) d(y)^{r-k}\,. \eeqano

Then,
$$ F(y) - P(x^*,y) = \sum_{j} \f_j(y) S_{j,0}(x^*,y)$$
and hence
$$ D^kF(y) - P_k(x^*,y) = \sum_{j} \sum_{i\leq k} \binomiale{k}{i} D^{k-i}\f_j(y) S_{j,i}(x^*,y)\,.$$

Therefore, choosing $\e$ sufficiently small: \beqano \|D^kF(y) -
P_k(x^*,y)\| &\leq& \sum_{j} \sum_{i\leq k} \binomiale{k}{i}
\|D^{k-i}\f_j(y)\|\cdot \|S_{j,i}(x^*,y)\|
\leq\\
&\leq& \sum_{j} \sum_{i\leq k} \binomiale{k}{i} M_{k-i}
d(y)^{-k+i}(C_k + \e) d(y)^{r-i} \leq\\
&\leq& E \,d(y)^{r-k} \leq \eta\, d(y)^{r-k-1}\,. \eeqano
\end{Proof}

\begin{Lem}\label{lemma5.4}
For every $\eta>0$, there exist $0<\d<1$ such that, for all $y\in
W_1\setminus A^*$, $x^*\in A^*$ and $x_0\in \dpr A^*$, we have
$$ \|P_k(x^*,y)-P_k(x_0,y)\|\leq \eta \|y-x_0\|^{r-k}\,, $$
whenever $k\leq r$ and
$$
\left\{ \begin{array}{l}
\|y-x^*\|<\a d(y)\\
\|y-x_0\|<\d\,.
\end{array}\right.
$$
\end{Lem}

\begin{Proof}
The proof goes as the one of lemma \ref{lemma5.2}, observing that
$\|x^*-x_0\|\leq (1+\a)\|y-x_0\|$ and
$$ P_k(x_0,y) - P_k(x^*,y) = \sum_{q\leq r-k} \frac{R_{k+q}(x^*,x_0)}{q!}(y-x)^q\,.$$
\end{Proof}

\begin{Claim}
For every $x_0\in \dpr A^*$ and $k=0,\,1,\,\ldots,\, r-1$:
$$D^kF(x_0)=f_k(x_0)\,.$$
Moreover, 
$D^{r-1}F$ is continuous at $x_0\in \dpr A^{*}$.
\end{Claim}

This claim follows easily from the lemmata above (see \cite{Abraham}, on page $128$, for more details).\\

This proves that $F\in C^{r-1}(\R^d)$ and completes the proof of i) and iv).\\
It remains to show that property vi) holds, namely that there exists
a constant $K>0$, such that $U(g(x))\leq K F(x)$ on $W_1$.
Obviously, this holds at every point in $A^*$, for every choice of
$K$ (since both functions vanish there).

\begin{Claim}
There exists a constant $K>0$, such that $\frac{U \circ g}{F} \leq
K$ on $W_1\setminus A^*$.
\end{Claim}

\begin{Proof}
Since $F>0$ on $W_1\setminus A^*$, it is sufficient to show that 
$\frac{U \circ g}{F}$ is uniformly bounded by a constant, as $d(y)$ goes to zero.

Let us start observing that, for $y\in K_j^{\l}$,
$$ \hat{P}_j(y) \geq C \|y-x_j\|^r \geq C d(y)^r\,;$$
therefore: \beqano F(y) &=& \sum_{j} \f_j(y) \hat{P}_j (y) \geq \\
&\geq& \sum_{j} \f_j(y) C d(y)^r = \\
&=& C d(y)^r\,.
 \eeqano

Moreover, if $x^* \in A^*$ such that $d(y)= \|y-x^*\|$, lemma
\ref{stimaresto} and \ref{lemma5.3} imply:
 \beqano
|U(g(y))-F(y)| &\leq& |U(g(y)) - P(x^*,y)| + |P(x^*,y)-F(y)| \leq \\
&\leq& C d(y)^r + Ed(y)^r = (C+E) d(y)^r\,. \eeqano

Hence,
 \beqano
\frac{U(g(y))}{F(y)} &=& \frac{U(g(y)) - F(y) + F(y)}{F(y)} \leq \\
&\leq& 1 + \frac{|U(g(y))-F(y)|}{F(y)} \leq \\
&\leq& 1 + \frac{(C+E) d(y)^r}{C d(y)^r} \leq \\
&\leq& 2 + \frac{E}{C} =: \tilde{K}\,.
 \eeqano
\end{Proof}

This proves property vi) and concludes the proof of the proposition.
\ep
\medbreak

\vspace{12pt}

\acks
I wish to thank John Mather for having introduced me to this area and suggested this problem.
I am very grateful to him and to Albert Fathi for their interest and
for several helpful discussions.\\

\nocite{*}
\bibliographystyle{plain}
\bibliography{biblioSorrentino}
{}

\end{document}